\documentstyle[12pt]{article}
\begin{document}
\input mssymb  
\newtheorem{dfn}{Definition}
\newtheorem{thm}[dfn]{Theorem}
\newtheorem{thrm}[dfn]{Theorem}
\newtheorem{chara}[dfn]{Characterization}
\newtheorem{mthm}[dfn]{Main Theorem}
\newtheorem{lem}[dfn]{Lemma}
\newtheorem{prop}[dfn]{Proposition}
\newtheorem{cor}[dfn]{Corollary}
\newtheorem{mlem}[dfn]{Main Lemma}
\newtheorem{sublem}[dfn]{Sublemma}
\newtheorem{ex}[dfn]{Example}
\newcommand{\fr}[1]{FillRad(#1)}
\newcommand{\sn}{S^{n}}
\newcommand{\slp}{S_{l}(p)}
\newcommand{\km}{$K_{M}$}
\newcommand{\exsists}{\exists}
\newcommand{\eps}{\varepsilon}
\renewcommand{\iff}{if and only if}
\renewcommand{\theequation}{\thethrm}
\title{The Radius Rigidity Theorem for Manifolds of Positive Curvature}
\author{Frederick Wilhelm
\thanks{Supported in part by a National Science Foundation
Postdoctoral Fellowship. }}
\date{January 19, 1995}
\maketitle

\begin{abstract}
Recall that the radius of a compact metric space $(X, dist)$ is given by
$rad\ X = \min_{x\in X} \max_{y\in X} dist(x,y)$.  In this paper we
generalize Berger's $\frac{1}{4}$-pinched rigidity theorem and
show that a closed, simply connected, Riemannian manifold with sectional
curvature $\geq 1$ and radius
$\geq \frac{\pi}{2}$ is either homeomorphic to the sphere or isometric to
a
compact rank one symmetric space.
\end{abstract}

The classical sphere theorem states that a complete, simply connected
Riemannian $n$-manifold with positive, strictly $1/4$-pinched  sectional
curvature is homeomorphic to $S^n$ ([Ber1], [K], and [Rch]).
   The weakly $1/4$-pinched case is covered by
\vspace{.1in} \newline
{\bf Berger's Rigidity Theorem}
 {\em  Let $M$ be a complete, simply connected
Riemannian $n$-manifold with sectional curvature, $1\leq sec\ M \leq 4$.
Then either
\begin{description}
\item[(\romannumeral1)]
	$M$ is homeomorphic to $S^n$, or
\item[(\romannumeral2)]
	$M$ is isometric to a compact rank one symmetric space.
\end{description}  }
([Ber2]) \newline

The hypotheses of Berger's Theorem imply (with
a lot of work) that the injectivity radius of $M$ satisfies
$inj\ M \geq \frac{\pi}{2}$ ([CG2] or [KS]).   The diameter therefore,
also
satisfies $diam\ M \geq \pi/2$, and the class of complete Riemannian
manifolds with
$$
sec \geq 1 \; \; \mbox{and} \; \; diam \geq \pi/2 \; \;
\hspace{2.2in} (*)
$$
 contains Berger's class.
The former class is in fact, much vaster, since it contains, for example,
metrics with arbitrarily small volume (see [Ber3] and Example 2.4 in
[GP1]).

On the other hand,
the set of smooth manifolds admitting metrics
satisfying $(*)$ is nearly the same as for Berger's class.
Indeed, in [GG1] Gromoll and Grove
extended Berger's Rigidity Theorem and the
Diameter Sphere Theorem ([GS]) in proving the \vspace{.1in} \newline
{\bf Diameter Rigidity Theorem:}
{\em Let $M$ be a complete, simply connected
Riemannian $n$-manifold with sectional curvature $sec\ M \geq 1$ and
diameter $diam\ M \geq \pi/2$.
Then either
\begin{description}
\item[(\romannumeral1)]
	$M$ is homeomorphic to $S^n$,
\item[(\romannumeral2)]
	$M$ is isometric to a compact rank one symmetric space, or
\item[(\romannumeral3)]
	$M$ has the cohomology algebra of the Cayley Plane, $CaP^2$.
\end{description} }

An open question regarding this
theorem is whether the possibility (\romannumeral3) can
be removed from the conclusion.
This seems to be a very difficult
problem; however, there is a natural hypothesis that falls between
those of the two rigidity theorems.  Observe that
the hypothesis $inj\ M \geq \pi/2$ (which is satisfied by Berger's class)
implies that given {\em any} point $x\in M$, there is a point $y\in M$
so that $dist (x,y) \geq \pi/2$.  This later condition can be expressed
succinctly in terms of a well known metric invariant called the
radius.
\begin{dfn}[Radius]\label{dfn of radius}
Let $(X,dist)$ be a compact metric space.  The {\em radius} of $X$ is
given by,
$$
rad\ X = \min_{x\in X} \max_{y\in X} dist(x,y).
$$
\end{dfn}
(The concept of radius was invented in [SY].  The name {\em radius}
was first used in [GP2].)

Clearly $inj\ M \geq \pi/2 \Rightarrow rad\ M \geq \pi/2 \Rightarrow
diam\ M \geq \pi/2$, suggesting the following generalization of
Berger's Rigidity Theorem. \vspace{.1in} \newline
{\bf Radius Rigidity Theorem:}
{\em Let $M$ be a complete, simply connected
Riemannian $n$-manifold with sectional curvature $sec\ M \geq 1$ and
radius $rad\ M \geq \pi/2$.
Then either
\begin{description}
\item[(\romannumeral1)]
	$M$ is homeomorphic to $S^n$, or
\item[(\romannumeral2)]
	$M$ is isometric to a compact rank one symmetric space.
\end{description} }


A crucial step in the proof of the Diameter Rigidity Theorem
is to show that if $M$ is not homeomorphic to $S^n$, then there are
certain points
$x$ whose unit tangent sphere
is mapped via $v \mapsto \exp_x \frac{\pi}{2}v$ onto the cut locus
of $x$, and that this map is a Riemannian submersion with connected
fibers.
Since the unit tangent sphere is isometric to the unit sphere $S^n \subset
{\Bbb R}^{n+1}$, the
classification theorem from [GG2] can be invoked.  It states that
up to isometric equivalence the only Riemannian submersions of
Euclidean spheres (with connected fibers) are the Hopf fibrations, except
possibly for
fibrations of the $15$-sphere by homotopy $7$-spheres.
It was shown in [GG1] that if the exception could be
removed from the submersion theorem in [GG2], then (\romannumeral3) can be
removed from the
statement of the Diameter Rigidity Theorem (see Remark 4.4 in [GG1]).
Although we have not been able to remove the exception from the submersion
classification, we have proved the following.
\begin{mlem}\label{1 point enough}
Let $S^{n}(r)$ denote  $\{ v\in {\Bbb R}^{n+1} \ | \ \|v\| = r \}$.
Let $\Pi : S^{15}(1) \longrightarrow V$ be a Riemannian submersion
with connected, $7$-dimensional fibers, and let $G$ be the set of
points $v\in V$ so that $\Pi^{-1}(v)$ is totally geodesic.

Then either $G$ is
discrete or $G$ is a totally geodesic
and isometrically embedded copy of $S^{l}(\frac{1}{2})$ for some
$1 \leq l \leq 8$.

Moreover, in case $G$ is discrete, either
\addtocounter{dfn}{1}
\begin{eqnarray} \label{rational distances}
dist( x,y) =  \frac{k\pi}{q}  \; \;  \mbox{for some $k, q \in {\Bbb N}$ so
that
					$q \not\equiv 0$ mod $4$, or}
					\nonumber     \\
	dist( x,y) =    \frac{\pi}{2}     \hspace{3.125in} \mbox{}
\end{eqnarray}
for all $x,y \in G$.
\end{mlem}

The Riemannian manifolds with
$$
sec\ M \geq 1, \; diam\ M\geq \frac{\pi}{2}, \;
\mbox{and nontrivial fundamental group} \; \; \; \; (**)
$$
 were completely classified in [GG1].
Naturally, the class with  $sec\ M \geq 1$, $rad\ M\geq \frac{\pi}{2}$,
and nontrivial fundamental group is contained in $(**)$.
It is not difficult to prove that this containment is proper.
\begin{thm}\label{pi_1}
Let $M$ be a closed, Riemannian $n$-manifold with sectional curvature
$sec\ M \geq 1$, radius $rad\ M \geq \frac{\pi}{2}$, and
nontrivial fundamental group $\Gamma$.  Then either
\begin{description}
\item[(\romannumeral1)]
	The universal cover $\tilde{M}$ of $M$ is isometric to
$S^{n}(1)$, and every orbit of the action of $\Gamma$ is contained in
a proper invariant totally geodesic subsphere, or
\item[(\romannumeral2)]
For some $d\geq 2$, $M$ is isometric to the ${\Bbb Z}_2$-quotient of
$CP^{2d-1}$
given by the involution
$$
[z_1, z_2, \ldots , z_{2d}] \mapsto [\bar{z}_{d+1}, \ldots ,
\bar{z}_{2d}, - \bar{z}_{1}, \ldots , - \bar{z}_d ]
$$
\end{description}
in homogeneous coordinates on $CP^{2d-1}$.

Moreover, all such spaces have $sec\ M \geq 1$ and $rad\ M \geq
\frac{\pi}{2}$.
\end{thm}

Recall ([W]) that a representation $\rho : \Gamma \longrightarrow O(n+1)$
is called fixed point free if and only if $S^n(1)/\rho(\Gamma)$ is
a space form.

The actions of the groups in (\romannumeral1) are necessarily reducible;
however, it is not immediately apparent (at least to the author)
exactly which (reducible) space forms satisfy the conclusion of
(\romannumeral1).
As a partial answer we will prove
\begin{thm}\label{big space forms}
Let $\rho: \Gamma \longrightarrow O(n+1)$ be a fixed point free
representation
that decomposes as a direct sum
$$
\rho_1 \oplus \rho_2 \oplus \cdots \oplus \rho_k
$$
of $k \geq 2$ irreducible representations.
\begin{description}
\item[(\romannumeral1)]
	A necessary condition for $S^n(1) / \rho(\Gamma)$ to
have radius $= \frac{\pi}{2}$ is that $\rho_i$ be equivalent to $\rho_j$
for some $i\not= j$.
\item[(\romannumeral2)]
	In case $\Gamma$ is abelian (\romannumeral1) is also a
sufficient condition.
\item[(\romannumeral3)]
	If $\Gamma$ is not abelian and $\tilde{\rho} : \Gamma
	\longrightarrow
O(d)$ is a fixed point free, irreducible representation, then
$rad\ \frac{S^{2d - 1}(1)}{ (\tilde{\rho} \oplus \tilde{\rho}) (\Gamma) }
< \frac{\pi}{2}$.
\item[(\romannumeral4)]
	If $rad\ S^n(1) / \rho(\Gamma) = \frac{\pi}{2}$ and $\sigma:
	\Gamma \longrightarrow
O(d)$ is another  fixed point
free representation of $\Gamma$, then $rad\ \frac{S^{n + d}(1) }{ (\rho
\oplus \sigma)
(\Gamma) } = \frac{\pi}{2}$.
\item[(\romannumeral5)]
	There is a $k_0$ (depending on $\Gamma$) so that if $k \geq k_0$,
then $rad\ S^n(1) / \rho(\Gamma) = \frac{\pi}{2}$.
\end{description}
\end{thm}

Given a smooth manifold $M$, the tangent and unit tangent bundles of $M$
will
be denoted by $TM$ and $SM$ respectively.  If $V\subset M$ is a smooth
submanifold, then the normal bundle of $V$ in $M$ will be denoted by
$NV$.
When there is no possibility of confusion we denote $S^n(1)$ by $S^n$.

For simplicity we abbreviate compact rank one symmetric space as CROSS.
All geodesics will be parametrized by arc length on $[0, \cdot]$ unless
otherwise indicated.

The remainder of the paper is divided into four sections and an appendix.
The first two sections contain the proof of the main lemma and a review of
certain
 material from [GG1].  The
Radius Rigidity Theorem is proved in section 3, and  Theorems \ref{pi_1}
and
\ref{big space forms} are proven in section 4.  In the appendix,
we give the proof of an inequality that is used in the proof of the main
lemma.
  \vspace{.1in}

\noindent
{\em Acknowledgements:}
It is my pleasure to thank
Carlos Duran,  Detlef Gromoll, Karsten Grove,  Luis Guijarro, and Vitali
Kapovich for stimulating conversations on this subject.  I especially
thank Karsten Grove for
discovering a mistake in an earlier draft and for offering other valuable
criticisms and
suggestions.  I would also like to thank Paula Bergen for her expert job
of copy editing.
Finally, I wish to acknowledge that the beautiful work in [GG1] and [GG2]
has been a true inspiration to me.

\section{Reflecting Good Points}
First we review some basic facts about Riemannian submersions.
Recall that if $\pi: M \longrightarrow B$ is a Riemannian submersion, then
the
vectors tangent to the fibers are called vertical vectors and the
vectors perpendicular to the fibers are called horizontal vectors.
We denote these two subbundles of $TM$ by $VM$ and $HM$ respectively.

The fundamental tensors of a submersion were defined in [O] as follows.
For arbitrary vector fields $E$ and $F$ on $M$ the tensor $T$ is defined
by
$$
T_E F = (\nabla_{E^v} F^v)^h + (\nabla_{E^v} F^h)^v,
$$
where the superscripts h and v denote the horizontal and vertical
parts of the vectors in question.  Note that the first summand
is the second fundamental form of a fiber applied to $E^v$ and
$F^v$, and the second term is the shape operator of a fiber applied to
$E^v$ and $F^h$.

The other fundamental tensor, $A$, is obtained by dualizing $T$, that is,
by
switching all horizontal and vertical parts in the definition of $T$.
Thus
$$
A_E F = (\nabla_{E^h} F^h)^v + (\nabla_{E^h} F^v)^h.
$$

	It is shown by O'Neill in [O], that all of the sectional
	curvatures of
$M$ can be written in terms
of  $A$, $T$, $\nabla A$, $\nabla T$, the sectional curvatures of $B$, and
the intrinsic
sectional curvatures of the fibers.  In particular, he proves that if $X$
and $Y$ are orthonormal horizontal
vector fields and $V$ is a unit vertical field, then
\addtocounter{dfn}{1}
\begin{description}
\item[Horizontal Curvature Equation]
\begin{equation}\label{horizontal curvature}
K(X,Y) = K( d\pi X, d\pi Y) - 3 \| A_X Y \|^2, \; \mbox{and}
\end{equation}
\addtocounter{dfn}{1}
\item[Vertizontal Curvature Equation]
\begin{equation}
K(X,V) = \langle (\nabla_X T)_V V, X \rangle + \| A_X V \|^2 - \| T_V
X\|^2.
\end{equation}
\end{description}

We refer the reader to [O] for the statements and proofs of the basic
facts
about $T$ and $A$ and other basic facts and definitions about Riemannian
submersions
that we will use freely and without further mention.

Now we begin the proof of the main lemma.
Let $\Pi$, $V$, and $G$ be as in Main Lemma \ref{1 point enough}.  We will
call the members of $G$ ``good points''.

\begin{lem}\label{wiedersehen}
\mbox{} \vspace{-.3in} \newline
\begin{description}
\item[ (\romannumeral 1) ]
 If $x \in G$, then there is a unique point $a(x) \in V$ at maximal
 distance
 from $x$, $ dist(x, a(x) ) = \frac{\pi}{2}$, and $a(x)$ is also in $G$.
\item[ (\romannumeral2) ]   $V$ is Wiedersehen at
$x$ and $a(x)$, i.e. the cut locus of $x$ is $a(x)$, and
the cut locus of $a(x)$ is $x$.  Furthermore, the fibers of
$\Pi$ are invariant under the antipodal map, and every geodesic in
$V$ is periodic with period $\pi$.
\end{description}
\end{lem}
{\em Remark:}
Gromoll and Grove have proven independently that the fibers of {\em any}
Riemannian submersion
with connected fibers of the $15$-sphere  are invariant under the
antipodal map
(even ones with $G= \emptyset$) ([GG3]).

\noindent
{\em Proof:}
First we review the notion of

\noindent
{\bf Holonomy Displacement}
	([H], p. 238, also [GG2], p. 150).
Let $\gamma$ be a geodesic in $V$.
If we consider all of the horizontal lifts of $\gamma$ to $S^{15}$,
then we obtain a family of diffeomorphisms, $\psi_{s,t} :
\Pi^{-1}(\gamma(s)) \longrightarrow \Pi^{-1}(\gamma(t) )$ given by
$\psi_{s,t}(z) = \gamma_z (t )$,
where $\gamma_z$ denotes the unique horizontal lift of $\gamma$
with $\gamma_z(s) = z$.

Now suppose that $F_x \equiv \Pi^{-1}(x)$ is totally geodesic.  Then
all horizontal geodesics emanating from $F_x$ are in a totally geodesic
$7$-sphere $F_{a(x)}$ at time $\pi/2$.  Hence $F_{a(x)}$ is also a fiber
of
$\Pi$, and $ \Pi (F_{a(x)}) $ is the desired point $a(x)$.
This proves (\romannumeral1).

Since every horizontal geodesic emanating from $F_x$ reaches $F_{a(x)}$ at
time
$\pi/2$, every geodesic emanating from $x$ reaches $a(x)$ at time
$\pi/2$, and hence is minimal up to time $\pi/2$.  Thus $V$ is Wiedersehen
at $x$ and by symmetry at $a(x)$.

It follows that reflection
in $x$ is a homeomorphism of $V$.  Hence reflection in
 $F_x$ is an
isometry of $S^{15}$ that maps fibers to fibers. Similarly, reflection
in $F_{a(x)}$ maps fibers to fibers.  But the composition of the two
reflections
is the antipodal map, $a$, of $S^{15}$.  So if we knew that
 the composition of
reflection in $x$ with reflection in $a(x)$ were the identity map of $V$,
then
we know that the fibers are invariant under the antipodal map.

 To establish this, let $r_x$, $r_{a(x)}$, $r_{F_x}$ and $r_{F_{a(x)} }$
 be the four reflections.
Note that $r_{F_x} \circ r_{F_{a(x)} } \circ r_{F_x} \circ r_{F_{a(x)} } =
a \circ a = id_{S^{15}}$.  So
$r_{x} \circ r_{a(x)} \circ r_{x} \circ r_{a(x)} = id_{V}$.
The differential $d (r_{x} \circ r_{a(x)} )_{a(x)}$ is therefore a linear
isometry of
$T_{a(x)}V$ whose square is the identity, and hence $T_{a(x)}V$ has a
basis of eigenvectors for  $d (r_{x} \circ r_{a(x)} )_{a(x)}$  with
corresponding eigenvalues of either $1$ or $-1$.
Suppose $v$ is an eigenvector whose eigenvalue is $-1$.  Then
$-v = d (r_{x} \circ r_{a(x)} )_{a(x)} v = d (r_{x})_{a(x)} -v$.
This implies that the reflection isometry $r_{x}$ fixes the
geodesic $c_{-v}: t \mapsto \exp_{a(x)} -tv$, which is absurd, since
$c_{-v} (\frac{\pi}{2}) = x$.   So the only possible
eigenvalue is $1$, and we can conclude that the fibers are indeed
invariant under the antipodal
map.  The invariance of the fibers under the antipodal map implies
immediately that
every geodesic in $V$ is periodic with period $\pi$. $\square$

We saw in the proof above that reflection in a totally geodesic fiber
is an isometry of $S^{15}$ that preserves the fibers of $\Pi$.
By using this fact over and over again, we can prove
\begin{lem}\label{reflection}
Let $x,a(x) \in V$ be good points at maximal distance. Let $z\in V
\backslash
\{ x, a(x)\}$ be another good point, and let $\gamma:[0, \infty)
\longrightarrow V$ be the unique geodesic  that passes through $x, z$, and
then $a(x)$ so that $\gamma(0) = x$.
Then $\gamma(k \cdot dist(x,z) )$ is a good point for all integers $k$.
In particular, if $dist(x,z)$ is an irrational multiple of $\pi$,
then all points along $\gamma$ are good.
\end{lem}

If $G$ has an accumulation point, then  using Lemmas \ref{wiedersehen} and
\ref{reflection}
and the fact that $G$ is closed, we see that $G$ contains
the image of an entire periodic geodesic of length $\pi$.
Thus in the indiscrete case it is enough to prove the following corollary
of the main lemma.

\begin{cor}\label{corollary of 1 point enough}
If $V^l \subset G_V$ is totally geodesic in $V$ and isometric
to $S^l(\frac{1}{2})$ for some $l \geq 1$, and if there is a good point
$x\in V \backslash V^l$, then there is
a totally geodesic set of good points, $V^{l+1} \subset V$, that contains
$V^l \cup \{ x\}$
and is isometric to
$S^{l+1}(\frac{1}{2})$.
\end{cor}
We will focus on the proof of (\ref{corollary of 1 point enough}) for
nearly all of the remainder of this section.

Note that $x \in G_V \backslash V^l$ implies $a(x) \in G_V \backslash
V^l$.
So by Lemma
\ref{wiedersehen}, $dist(x, \cdot)$ is smooth along $V^l$.  If
 $dist(x, \cdot)|_{V^l}$ is not constant,
there is therefore some open set $O\subset V^l$ for which the subset
$I_O = \{ q\in O \ | \ dist(x,q) \; \mbox{is an irrational multiple
of $\pi$ } \}$ is a dense $G_{\delta}$.  By Lemma \ref{reflection},
the geodesics between $x$ and the points in $I(O)$
consist of good points.  By continuity, then,
all geodesics between $x$ and points of $O$ consist of good points.
 Since $x$ is a good point and $a(x)\not\in V^l$, there is a unique
 minimal geodesic between $x$
and every point in $O$.   Let $C(O)$ denote the union of these geodesics.
Then $C(O) \backslash \{ x, a(x) \}$ is a smooth, $(l +1 )$-dimensional
submanifold
of $V$ composed entirely of good points.

Now consider a point $y \in O$.  Let $V^{l+1}$ be the
image of the set of geodesics emanating from
$y$ which are initially tangent to $C(0)$.  It follows from
Lemma \ref{reflection}
that $V^{l+1}$ consists entirely of good points, and by
Lemma \ref{wiedersehen} it is homeomorphic to $S^{l+1}$.
To help understand the infinitesimal geometry of $V^{l+1}$
we prove the following.

\begin{lem}\label{constant curvature}
If $\gamma: [0, \pi] \longrightarrow V$ is a geodesic whose image consists
entirely
of good points, then all radial sectional curvatures of $V$ along $\gamma$
are constant and equal to 4.
\end{lem}
{\em Proof:}
Let $\tilde{\gamma}$ be a horizontal lift of $\gamma$ to $S^{15}$.
Let $X$ and $\tilde{X}$ denote the tangent fields to
$\gamma$ and $\tilde{\gamma}$ respectively, and let $V$ be a
vertical unit field along $\tilde{\gamma}$ so that $(\nabla_{\tilde{X}}
V)^v =0$.  Then, using the equation for the
vertizontal curvatures, we find that $K(\tilde{X}, V)$ along
$\tilde{\gamma}$ is,
\addtocounter{dfn}{1}
\begin{eqnarray}\label{vertizontal}
1 \equiv K(\tilde{X}, V) = \langle (\nabla_{\tilde{X}} T)_V V, \tilde{X}
\rangle + \| A_{\tilde{X}} V \|^2
- \| T_V \tilde{X} \|^2 = \| A_{ \tilde{X} } V\|^2.
\end{eqnarray}
It follows from (\ref{vertizontal}) that
\addtocounter{dfn}{1}
\begin{eqnarray}\label{bijective}
\mbox{the map $v \mapsto A_{\tilde{X}} v$ from
$VS^{15}$ to $HS^{15} \cap \tilde{X}^{\perp}$ is bijective.}
\end{eqnarray}
Combining this and
(\ref{vertizontal}) we can show that
\addtocounter{dfn}{1}
\begin{eqnarray}\label{horizontal A-tensor value}
\| A_{\tilde{X}} y\| \equiv 1 \; \; \mbox{for all unit vectors $y \in
HS^{15} \cap
\tilde{X}^{\perp}$.}
\end{eqnarray}
Indeed, if $\langle A_{\tilde{X}} y, v\rangle$ were bigger than $1$ for
some unit vector
$v\in VS^{15}$, then we would have $| \langle A_{\tilde{X}} v, y \rangle |
=
| \langle v, A_{\tilde{X}} y \rangle | > 1$, contrary to
(\ref{vertizontal}).
On the other hand, by (\ref{vertizontal}) and (\ref{bijective}),
 $y = A_{\tilde{X}} v$ for some unit vector $v \in VS^{15}$,
thus $| \langle A_{\tilde{X}} y, v \rangle | = | \langle y, A_{\tilde{X}}
v \rangle | = 1$, and
$\| A_{\tilde{X}} y \| \geq 1$ as well.

Let $Y$ be a unit field along $\gamma$ that is perpendicular to
$X$, and let $\tilde{Y}$ denote the horizontal lift of $Y$.
Then it follows from (\ref{horizontal A-tensor value}) that
$$
K(X, Y) = K(\tilde{X},\tilde{Y}) + 3 \| A_{\tilde{X}} \tilde{Y} \|^2 =
1 + 3* 1 = 4.
$$
$\square$

The proof of Corollary \ref{corollary of 1 point enough} in
case  $dist(x, \cdot)|_{V^l}$ is not constant is
completed by applying the following result to $V^{l+1}$.
\begin{lem}\label{V^{l+1} is totally geodesic}
Let $P\subset T_yV$ be a $k$-dimensional subspace of a tangent fiber of
$V$
such that $W\equiv \exp_y P$ is a subset of $G_V$.
Then $W$ is totally geodesic and isometric to $S^k(\frac{1}{2})$.
\end{lem}
{\em Proof:}
Let
$\iota: T_y V \longrightarrow T_zS^{8}(\frac{1}{2})$ be a linear
isometry.
It follows from Lemma \ref{constant curvature}
that the differentials of $\exp_z \circ \iota \circ
\exp^{-1}_{y}$ are isometries along $V^{l+1}$, and it follows from Lemma
\ref{wiedersehen} that $\exp_z \circ \iota \circ \exp_{y}^{-1}|_V$
is an embedding whose image is isometric to $S^{k}(\frac{1}{2})$.
$\square$

To prove Corollary \ref{corollary of 1 point enough} it
remains to consider the case when the restriction of
$dist(x, \cdot )$ to $V^{l}$ is constant.  In this case
it turns out that
\addtocounter{dfn}{1}
\begin{eqnarray}\label{constant distance}
dist(x,\cdot )|_{V^{l}} \equiv \pi/4.
\end{eqnarray}
  To see this, first note that
for $z\in V^l$, $\pi/2 = dist(z,a(z) )\leq
dist(z, x) + dist( x,a(z) ) =
2dist( z , x )$.
On the other hand $d_{Haus} (S^{15}, \pi^{-1} (V^l) )
\leq \pi/4$, so $dist( x, V^l ) \leq \pi/4$.

By combining (\ref{constant distance}) with Lemma \ref{wiedersehen},
we see that all geodesics from $x$ that pass through $V^l$
go through good points at times $0$, $\frac{\pi}{4}$, $\frac{\pi}{2}$,
and $\frac{3\pi}{4}$.   For the rest of this section we will study
geodesics with this property.  Our goal will be to prove that all of the
points along such a geodesic are good.  The first step is to estimate
the average Ricci curvature along $\gamma$. \vspace{.1in} \newline
\begin{lem}\label{horizontal average lemma}
Let $\gamma :[0, \pi] \longrightarrow V$ be a geodesic so that the points
$\gamma(0)$, $\gamma(\pi/4)$,
$\gamma(\pi/2)$, and $\gamma (3\pi/4 )$ are good, and let $\{ E_i
\}_{i=1}^{7}$
be an orthonormal collection of parallel unit normal fields along
$\gamma$.  Then
\addtocounter{dfn}{1}
\begin{equation}\label{average base bound}
\int_{0}^{\pi}  Ric (\dot{\gamma}, \dot{\gamma} )   \leq 28 \pi, \;
\mbox{and}
\end{equation}
\addtocounter{dfn}{1}
\begin{equation}\label{average A-tensor bound}
\int_{0}^{\pi} \Sigma_{i=1}^{7} \| A_{  \dot{ \tilde{\gamma}} }
\tilde{E}_i \|^2  \leq 7 \pi,
\end{equation}
where $\dot{ \tilde{\gamma} }$ and $\tilde{E}_i$ denote the horizontal
lifts
of $\dot{\gamma}$ and $E_i$ respectively.
\end{lem}
{\em Proof:}
Inequality (\ref{average A-tensor bound}) is clearly a consequence of
(\ref{average base bound}),
the equation for horizontal curvatures, and the fact the the curvature of
$S^{15}$ is identically $1$.  So it suffices to prove (\ref{average base
bound}).

By Lemma \ref{wiedersehen}
there are no conjugate points along $\gamma$ prior to time $\pi/2$.  Using
this and a
``Bonnet-Meyers" type of argument we can show,
\addtocounter{dfn}{1}
\begin{eqnarray}\label{component average}
 \int_{0}^{\pi/2} K(\dot{\gamma}, E_i) \sin^2 2t \leq 4 \int_{0}^{\pi/2}
 \sin^2 2t
\end{eqnarray}
for all $i = 1, \ldots, 7$.

Suppose (\ref{component average}) is false.  Then for some $l <
\frac{\pi}{2}$
and $i= 1, \ldots, 7$,
\addtocounter{dfn}{1}
\begin{equation}\label{bad average}
 \int_{0}^{l} \sin^2 \frac{\pi t}{l} K(\dot{\gamma}, E_i) >
\frac{\pi^2}{l^2} \int_{0}^{l} \sin^2 \frac{\pi t}{l}.
\end{equation}

Now set $W_i = \sin(\frac{\pi t}{l}) E_i(t)$, and compute the index:
\begin{eqnarray*}
I(W_i,W_i) = \int_{0}^{l} \langle \nabla_{\dot{\gamma}}
\sin(\frac{\pi t}{l}) E_i,  \nabla_{\dot{\gamma}} \sin(\frac{\pi t}{l})
E_i
\rangle - \sin^2 (\frac{\pi t}{l}) \langle R( E_i,
\dot{\gamma})\dot{\gamma},
 E_i\rangle dt= \\
 \int_{0}^{l} - \langle
\sin(\frac{\pi t}{l}) E_i, \nabla_{\dot{\gamma}}  \nabla_{\dot{\gamma}}
\sin(\frac{\pi t}{l}) E_i
\rangle - \sin^2 (\frac{\pi t}{l}) \langle R( E_i,
\dot{\gamma})\dot{\gamma},
 E_i \rangle dt = \\
\int_{0}^{l}  \frac{\pi^2}{l^2} \sin^2(\frac{\pi t}{l}) \langle
 E_i,  E_i
\rangle - \sin^2 (\frac{\pi t}{l}) \langle R( E_i,
\dot{\gamma})\dot{\gamma},
 E_i \rangle dt = \\
\int_{0}^{l}  \sin^2(\frac{\pi t}{l})( \frac{\pi^2}{l^2}
 -  \langle R( E_i, \dot{\gamma})\dot{\gamma},
 E_i \rangle ) dt
\end{eqnarray*}
If (\ref{bad average}) holds for some $i= 1, \ldots 7$, it follows that
$I(W_i,W_i) < 0$, implying that $\gamma|_{[0,l]}$ is not minimal, a
contradiction.

 Applying the same argument to $\gamma|_{[\pi/4, 3\pi/4]}$,
 $\gamma|_{[\pi/2, \pi]}$,
and $\gamma|_{[3\pi/4, 5\pi/4]}$ shows
$$
 \int_{\pi/4}^{3\pi/4} K(\dot{\gamma}, E_i) \cos^2 2t \leq 4
 \int_{\pi/4}^{3 \pi/ 4} \cos^2 2t,
$$
$$
 \int_{\pi/2}^{\pi} K(\dot{\gamma}, E_i) \sin^2 2t \leq 4
 \int_{\pi/2}^{\pi} \sin^2 2t, \; \; \mbox{and}
$$
$$
 \int_{3\pi/4}^{5\pi/4} K(\dot{\gamma}, E_i) \cos^2 2t \leq 4
 \int_{3\pi/4}^{5\pi/4} \cos^2 2t
$$
for all $i=1, \ldots 7$.
Combining these with (\ref{component average}) and using the fact that
$\gamma$ is periodic with period $\pi$ yields (\ref{average base bound}).
								$\square$

   From (\ref{average A-tensor bound}) we get \vspace{.1in} \newline
\begin{lem}\label{appendix lemma}
Let $\gamma$ be as in Lemma \ref{horizontal average lemma}.
Let $\{v_i \}_{i=1}^7$ be an orthonormal basis for $VS^{15}_{\gamma(0)}$,
and
let $V_i$ be an extension of $v_i$ to a vertical field such that
$(\nabla_{\dot{\tilde{ \gamma } } }  V_i)^v =0$.
Then
\addtocounter{dfn}{1}
\begin{eqnarray}\label{average vertical A-tensor}
\int_{\tilde{\gamma}} \Sigma_{i=1}^{7}  \| A_{ \dot{ \tilde{\gamma} } }V_i
\|^2
\leq 7\pi.
\end{eqnarray}
\end{lem}
The verification of Lemma \ref{appendix lemma} is a lengthly but
rather routine exercise in linear
algebra, so we defer it to the appendix.

The proof of Corollary \ref{corollary of 1 point enough}
is completed by combining Lemmas \ref{constant curvature} and
\ref{V^{l+1} is totally geodesic} and Equation \ref{constant distance}
 with the following result.
\begin{lem}\label{four good points}
If $\gamma :[0, \pi] \longrightarrow
V$ is a geodesic so that the fibers $p^{-1}(\gamma(0))$,
$p^{-1}(\gamma(\pi/4))$, $p^{-1}(\gamma(\pi/2) )$, and
$p^{-1}( \gamma (3\pi/4 ) )$ are totally geodesic, then all of the fibers
$p^{-1}(\gamma(t))$ are totally geodesic.
\end{lem}
{\em Proof:}
Let $\{ V_i\}_{i=1}^{7}$ be as in the statement of (\ref{appendix
lemma}).
Averaging vertizontal curvatures along $\tilde{\gamma}$ we find
\addtocounter{dfn}{1}
\begin{eqnarray}\label{vertizontal average}
\int_{0}^{\pi} \Sigma_{i=1}^{7} \langle R(V_i, \dot{\tilde{\gamma} } ),
\dot{\tilde{\gamma} } ,V_i
\rangle = \int_{0}^{\pi} \Sigma_{i=1}^{7} \dot{\tilde{\gamma} } \langle
T_{V_i} V_i, \dot{\tilde{\gamma} }  \rangle + \| A_{\dot{\tilde{\gamma} }
} V_i \|^2 - \| T_{V_i} \dot{\tilde{\gamma} } \|^2 dt
\end{eqnarray}
The first term on the right is equal to
\begin{eqnarray*}
\Sigma_{i=1}^{7} \langle T_{V_i} V_i, \dot{\tilde{\gamma} }
\rangle|^{\pi}_{0}.
\end{eqnarray*}
All of these terms are zero, since $\gamma(0) = \gamma(\pi)$ is a good
point.
Thus (\ref{vertizontal average}) becomes
\begin{eqnarray*}
7 \pi
= \int_{\gamma} \Sigma_{i=1}^{7} \langle R(V_i, \dot{\tilde{\gamma} } ),
\dot{\tilde{\gamma} } ,V_i
\rangle = \int_{\gamma} \Sigma_{i=1}^{7}  \| A_{\dot{\tilde{\gamma} } }
V_i \|^2 - \| T_{V_i} \dot{\tilde{\gamma} } \|^2 dt.
\end{eqnarray*}
Combining this with (\ref{average vertical A-tensor}) shows
\begin{eqnarray*}
 \int_{\gamma} \Sigma_{i=1}^{7}  \| T_{V_i} \dot{\tilde{\gamma} } \|^2 dt
 = 0.
\end{eqnarray*}
Thus $T_{V_i} \dot{\tilde{\gamma} } \equiv 0$, and hence
\addtocounter{dfn}{1}
\begin{equation}\label{partially geodesic}
T_v \dot{\tilde{\gamma} }  = 0 \;
\mbox{for all $v\in VS^{15}|_{ \dot{\tilde{\gamma} } }$.}
\end{equation}

 It turns out that the condition
(\ref{partially geodesic}) implies that all of the
Holonomy Displacement maps for $\gamma$
are isometries.  This and the fact that the fiber $\pi^{-1} (\gamma(0))$
is totally geodesic yields the conclusion of Lemma \ref{four good
points}.

So it remains to see that the above Holonomy Displacement maps,
$\psi_{s,t}$,
are isometries.  Consider a curve $c:[0,l] \longrightarrow \Pi^{-1}
(\gamma(t_0) )$
with $\| \dot{c}\| \equiv 1$ and the variation $W(s,t)$ of $c$ that is
given by $W(s,t) =
\psi_{t_0, t} (c(s) )$.  The variation field of $W$ along $c$ is
the horizontal lift of $\dot{\gamma}(t_0)$.  Denote it by $\tilde{X}$.
By the first variation formula, we have
\begin{eqnarray*}
\frac{d}{dt} Length [ W( \cdot, t) ] =
\langle \tilde{X}, \dot{c} \rangle|^{l}_{0} - \int^{l}_{0} \langle
\tilde{X},
\nabla_{\dot{c}} \dot{c} \rangle = \\
\int^{l}_{0} \langle T_{\dot{c}} \tilde{X},
 \dot{c} \rangle = 0.
\end{eqnarray*}
The second equality is due to the properties of $T$ and the facts that
$\tilde{X}$ is horizontal and $\dot{c}$ is vertical.  The last equality
follows from (\ref{partially geodesic}).   $\square$

To complete the proof of the main lemma, we note that
if $G$ is discrete, then by Lemmas \ref{wiedersehen}, \ref{reflection},
and \ref{four good points}, the equations in (\ref{rational distances})
hold.

\section{Review of the Diameter Rigidity Theorem}

If $M$ satisfies the hypotheses of the Radius Rigidity Theorem,
then $M$ also satisfies the hypotheses of the Diameter Rigidity Theorem,
so the only way $M$ can fail to satisfy the conclusion of the
Radius Rigidity Theorem is if it has the cohomology algebra of
$CaP^2$.  We assume throughout sections 2 and 3 that $sec\ M \geq 1$,
\label{exceptional}
$Rad\ M \geq \pi/2$, $\pi_1(M) = \{ e\}$, and $H^*(M) \cong H^*(CaP^2)$,
and we attempt
to show that $M$ is isometric to $CaP^2$.

By the Diameter Sphere Theorem ([GS]), $diam\ M = \pi/2$.
We would like to focus on this property for awhile; so let $N$ be a
Riemannian $n$-manifold with $sec\ N \geq 1$, $\pi_1(N) \cong \{ e\}$, and
$diam\ N =\pi/2$,
that is not homeomorphic to $S^n$.
Many basic aspects of the geometry of $N$ can be described in
terms of so called {\em dual sets} ([GG1]).  (Cf also [Sa], [Sh], and
[SS].)
\begin{dfn}[Dual Sets]\label{dual sets}
For any subset $B \subset N$, the dual set of $B$ is,
$$
B' = \{ x\in N \ | \ dist(x, B) = \pi/2 \}
$$
\end{dfn}
The following properties of dual sets were observed in [GG1] (cf also
 [Sa], [Sh], and [SS]).
\begin{description}
\item[(\romannumeral1)]
	$B'$ is totally $\pi$-convex,  that is, any geodesic of length
	strictly less than
$\pi$ whose end points lie in $B'$ lies entirely in $B'$.
\item[(\romannumeral2)]
	 $B \subset B''$.
\item[(\romannumeral3)]
	If $A \subset B$, then $A' \supset B'$.
\item[(\romannumeral4)]
	$B' = B'''$.
\end{description}

It follows from (\romannumeral1) and [CG1] that $B$ is a topological
manifold with (possibly empty)
boundary and smooth, totally geodesic interior.
If we start with a set $B$ so that $B' \not= \emptyset$ and set $A= B'$,
then
$A = (A')'$.  Thus
$$
A= \{ x\in N \ | \ dist(x, A') = \pi/2 \} \; \mbox{and} \;
A' = \{ x\in N \ | \ dist(x, A) = \pi/2 \}, \; \mbox{and}
$$
$A$ and $A'$ are called a {\em dual pair}.

The proof in [GG1] proceeds from this point to use comparison theory
and other geometric and topological tools  to
argue that the geometry of $N$ is more and more like the geometry of
a CROSS.  For example, it is shown that $\partial A = \partial A'=
\emptyset$,
that cutlocus$(A)= A'$ and cutlocus$(A') = A$, and that for any $p \in A$
the
map $\Pi_{p} : UNA_{p} \longrightarrow A'$ from the unit normal sphere
to $A$ at $p$ to $A'$ given by $\Pi_{p}(u) = \exp (\frac{\pi}{2} u)$
is a Riemannian submersion with connected fibers.  This allows them to
apply the classification
theorem in [GG2] and conclude that $\Pi_{p}$ is isometrically
equivalent to a Hopf fibration (except possibly if the fibers are
$7$-dimensional).
The proof is completed with further comparison arguments.  The exception
to
the conclusion is accounted
for by the fact that the classification in [GG2] is not quite
complete.  It leaves open the
possibility of nonstandard Riemannian submersions of the $15$-sphere by
homotopy $7$-spheres.  On the other hand, using arguments from
[GG1] it is easy to prove that this is the only possible obstruction.
\begin{prop}\label{1 nice fibration}
If $N$ has a dual pair $(A, A')$ such that one of the submersions
$\Pi_{p}$ is isometrically
equivalent to a Hopf fibration, then $N$ is isometric to a CROSS.
\end{prop}
{\em Proof:}
Say $p \in A$,  and $\Pi_{p}$ is isometrically
equivalent to a Hopf fibration.   It was shown in [GG1] (p. 236) that it
is enough to find a
dual pair $\{ q \}$, $\{ q \}'$, where $\{ q\}$ is a singleton and $\Pi_q$
is isometrically equivalent to a Hopf fibration.  So we may assume that $A
\not= \{ p\}$.
By the Diameter Rigidity Theorem, we may assume that $N$ has the
properties of the
possibly exceptional manifold $M$, on page \pageref{exceptional}, $sec\ N
\geq 1$,
$Rad\ N \geq \pi/2$, $\pi_1(N) = \{ e\}$, and $H^*(N) \cong H^*(CaP^2)$.
So we can refer to $N$ as $M$.
We also know that $A'$ is isometric to a CROSS, $P^{m}(K)$.
It was observed in [GG1] (p. 236) that the dual set $B$ (in $A'$) of any
singleton $\{ x \} \subset A'$  is isometric to $P^{m-1}(K)$, and that
the double dual of $\{x\}$ (in $M$) is again $\{ x\}$.
It follows from the convexity properties of $A'$ that the fibers of the
submersion
$SA_{x}' \longrightarrow B$ are also fibers of the submersion $\Pi_x :
SM_x \longrightarrow
\{ x\}'$, and it follows from our simplifying assumptions that the
dimension of these fibers
is $< 7$.  Therefore the submersion $\Pi_x$ is equivalent to a Hopf
fibration and
$M$ is isometric to a CROSS. $\square$

We now restrict our attention to the possibly exceptional manifold $M$.
\begin{prop}\label{covering by dual pairs}
\mbox{} \vspace{-.3in} \newline
\begin{description}
\item[(\romannumeral1)]
	The set of dual pairs is a covering of $M$.
\item[(\romannumeral2)]
	Every dual pair consists of a singleton and a set that is
	homeomorphic to
$S^8$.
\item[(\romannumeral3)]
If $(p, V)$ and $(q,W)$ are distinct dual pairs, then $V \cap W $ is a
point.
\end{description}
\end{prop}
{\em Proof:}
	(\romannumeral1) is an immediate consequence of properties
	(\romannumeral2) and
(\romannumeral4) of dual sets and the fact that $rad\ M = \pi/2$.

To prove (\romannumeral2) first note that if $(p, V)$ is a dual pair and
$V$ is not $8$-dimensional,
then the Riemannian submersion $\Pi_{p}: SM_p \longrightarrow V$ is
isometrically equivalent to
a Hopf fibration; so $M$ is isometric to a CROSS.   If $V$ is
$8$-dimensional, then the fibers
of $\Pi_{p}$ are homotopy $7$-spheres (see Theorem 5.1 in [Br]).  It
follows from the long exact
homotopy sequence of the fibration $\Pi_p$ that $V$ is a homotopy
$8$-sphere, and hence a topological $8$-sphere.  Finally, if there is a
dual pair $(A, A')$ so that neither $A$ nor $A'$
is a point, then $1 \leq dim\ NS_{p} \leq 14$, and the submersion $\Pi_p$
is equivalent to a Hopf
fibration.

To prove (\romannumeral3) observe that
since $sec\ M> 0$ and $dim\ V + dim\ W = dim\ M$, a Synge Theorem type of
argument shows that $V\cap W\not= \emptyset$  (see [F] and also
Proposition 1.4 in [GG1]).
Next observe that
$V\cap W = \{ p, q\}' $, so $(V\cap W, (V\cap W)' )$ is a dual pair.  By
(\romannumeral2),
 one of these dual sets is a point.  Since $p, q \in (V\cap W)' $, we
 conclude that $(V\cap W)$ is
a point.  $\square$

If $(p, V)$ is a dual pair, then we will (optimistically) refer to $V$
as a {\em Cayley line}.  This name is partially justified by the fact that
once we have proven that $M$ is isometric to $CaP^2$ we will know that
all of these $V$'s are isometric to $CaP^1$.

\section{Intersecting Cayley Lines}
In this section we prove the Radius Rigidity Theorem.

	If $(p, V)$ is a dual pair, then we have seen that it is enough
to show that the submersion $\Pi_p :S_p \longrightarrow V$ is
isometrically equivalent to the Hopf fibration $S^7 \hookrightarrow
S^{15} \longrightarrow S^8$.  This holds if
its fibers are totally geodesic (see [Rj]).  Roughly speaking, the
strategy of our
proof is to find dual pairs $(p, V)$ so that $\Pi_{p}$ contains
more and more totally geodesic fibers.  Our method for finding
totally geodesic fibers will be to find more and more ``good points'' in
$M$.
\begin{dfn}[Good Point] \label{good point}
If $(p, V)$ is a dual pair, then we shall call a point $x \in V$ good
if and only if $\Pi_{p}^{-1}(x)$ is totally geodesic.
 A point $m \in M$
will be called good if and only if $m\in G_W$ for some Cayley line
$W \subset M$.  The sets of good points in $V$ and $M$ will be denoted by
$G_V$ and $G_M$
respectively.
\end{dfn}
The fact that $G_M$ is rather large is a consequence of Proposition
\ref{covering by dual pairs} and the next result.
\begin{prop}\label{properties of good points}
Let $(p, V)$ be a dual pair. \vspace{-.3in} \newline
\begin{description}
\item[(\romannumeral1)]
	A point $x\in V$ is good if and only if there is a
Cayley line $W$ so that $V\cap W = \{ x\}$.
\item[(\romannumeral2)]
	$G_M$ is closed.  In fact, $m \in G_M$ if and only if there are
	points $x, y \in M$
so that $dist(x,m) = dist(m,y) = dist(x,y) = \frac{\pi}{2}$.
\end{description}
\end{prop}
{\em Remark:} Gromoll and Grove proved independently that
 $G_M \not= \emptyset$ ([GG3]). \newline
{\em Proof:}
If there is a dual pair $(z, W)$ so that $W\cap V = \{ x\}$,
then $p$ and $z$ are distinct
points in $\{ x\}'$, so $\{ x\}'$ is a Cayley line
and hence intersects $V$ in a single point $y$.
Since $p$ and $x$ are distinct points of $\{ y\}'$,
$\{ y\}'$ is a Cayley line.  It follows that the set
of minimal geodesics from $p$ to $x$ is contained in $\{ y\}'$.
Thus $\Pi^{-1}_{p}(x)$ is contained in the unit tangent sphere $S\{
y\}_{p}'$
to $\{ y\}'$ at $p$.  But since both of these sets are homotopy
$7$-spheres
they must coincide.
compact
Since $S\{ y\}_{p}'$ is totally geodesic in $SM_p$, $\Pi^{-1}_{p}(x)$ is
as
well.  This proves the ``if'' part of (\romannumeral1).

On the other hand,
if $x\in G_V$, then by Lemma \ref{wiedersehen} there is a unique point
$a(x) \in V$ so that $dist(x, a(x) ) =\frac{\pi}{2}$.  Since
$x, p \in \{ a(x) \}'$, $\{ a(x) \}'$ is a Cayley line.
By Proposition \ref{covering by dual pairs}, $x = \{a(x) \}' \cap V$.
This proves the ``only if'' part of (\romannumeral1).
Since $dist(x, a(x) ) = dist (a(x), p) = dist (p,x) = \frac{\pi}{2}$ it
also proves
the ``only if'' part of (\romannumeral2).

To prove the ``if'' part of (\romannumeral2) note that $x,y \in \{ m\}'$,
$m,y \in \{x\}'$, and $m,x \in \{y\}'$.  So $\{ m\}'$, $\{ x\}'$, and
$\{ y\}'$ are all Cayley lines, and $x$, for example, is good since
$x = \{ m\}' \cap \{ y\}'$.
     $\square$

The Radius Rigidity Theorem would follow if we could show that
there is a Cayley line
$V$ so that every point in $V$ is good.  We will do this by
finding Cayley lines with
good points in sets that are isometric to spheres of
constant curvature $4$ of progressively higher and higher dimension.
Since each point in $M$ lies on at least one Cayley line, we can certainly
find a countably infinite family of Cayley lines $\{ V_i
\}_{i=1}^{\infty}$.
Next we observe that there is a Cayley line $W$ so that
$$
card\ \{ W \cap V_i \}_{i=1}^{\infty}
$$
is infinite.  Indeed if $card\ \{ W \cap V_i \}_{i= 1}^{\infty}$ were
finite,
then there would be an infinite set $\{ V_{i_j}\}_{j=1}^{\infty}$
and a point $x$ so that $W\cap V_{i_j} = \{ x\}$ for all $j$.  But then
$\{ x\}'$ is a Cayley line and the points
$\{ x\}' \cap V_{i_j}$ must all be distinct. So we can find a Cayley line
(let's call it $V$) with infinitely many distinct good points.
It follows that the set of good points in $V$ contains an accumulation
point
and hence, using Lemma \ref{reflection} and the fact that $G_V$ is closed,
 the image of an entire geodesic.

To prove the Radius Rigidity Theorem
 we argue by contradiction.  It follows from the main lemma
that the set of good points in each Cayley line is either discrete,
an entire geodesic, or a sphere of constant curvature $4$.  Let $V$ be
a Cayley line whose set of good points has maximal dimension, $d$.
We've seen that $d \geq 1$, and if the Radius Rigidity Theorem were
false, then we would know $d \leq 7$.  Consider the configuration $C$
consisting
of all Cayley lines of the following types:
\label{configuration}
\begin{description}
\item[type 1]
	$V$,
\item[type 2]
	All the Cayley lines between the good points of $V$ and $V'$,
\item[type 3]
	For each $W$ of type 2 we also include all of the Cayley lines
between each of the good points of $W$ and $W'$ that are neither of type 1
nor of type 2 .
\end{description}
We point out that if $W$ is of type 2, then $W'$ is a good point of
$V$, and if $U$ is a line of type 3 between a good point of $W$ and
$W'$, then $U'$ is a good point of $W$.

 Suppose we could find a Cayley line $Z$ that is not included
in the configuration above.  Then either,
\begin{eqnarray*}
Z \cap V \not\in G_V   \;    \mbox{or}     \\
Z\cap V \in G_V
\end{eqnarray*}
But neither of these is possible.  The first can not occur because $G_V$
consists of all
of the good points of $V$.  On the other hand, if $Z\cap V \in G_V$, then
$(Z\cap V)'$ is a line of type 2, and $Z$ is a Cayley line
between $(Z\cap V)'$ and $Z\cap V$, implying that $Z$ is of
type 3, a contradiction.  Therefore,
\addtocounter{dfn}{1}
\begin{eqnarray}\label{C is everything}
\mbox{$C$ contains all of the lines of $M$.}
\end{eqnarray}

Next we prove
 \begin{lem}\label{d is constant}
We may assume that $dim\ G_U = d$ for every line $U$ in the
configuration.
\end{lem}
{\em Proof:}
We prove this in a step by step manner.

We know that there is at least one line of type 3, since otherwise the
configuration would only be $8$-dimensional in a neighborhood of a bad
point
of $V$ and hence would not cover $M$.  Since all of the
lines of type 2 intersect at $V'$, they must intersect at distinct
points of each line of type 3.  Thus $dim\ G_U \geq 1$ for all lines of
type 3, and
if $U$ is a line of type 3, then there are infinitely many lines between
$U$ and $U' \in W_0$, where $W_0$ is a line of type 2.  Since all of these
lines
intersect at $U'$ they must intersect each line of type 2 (other than
$W_0$) in
infinitely many places.  Therefore $dim\ G_W \geq 1$ for all lines $W\not=
W_0$
of type 2.  Since the set of all good points in $M$ is closed,
 $dim\ G_{W_0} \geq 1$ as well.

For each point
$v \in G_V$, the set $L_v \equiv
\{\mbox{lines $U$ in the configuration} \ | \ v\in U \}$
can be topologized by declaring that it is homeomorphic to
$G_{ \{ v\}' }$.
We will show that for each $v\in G_V$, $L_{v}^d \equiv \{ \mbox{lines $U$
in }
L_v \ | \  dim\ G_U = d \}$ is both closed and open.  Since $V \in
L_{v}^{d}$
and $\cup_{v\in G_v} L_v = C$,
it will follow that $dim\ G_U = d$ for every line in the configuration.
Let $\{ U_i\}$ be a sequence in $L_{v}^{d}$ which is converging to a line
$U$ in $L_v$.  Then by passing to a subsequence if necessary, we may
assume that
$\{G_{U_i} \}$ converges (in the classical Hausdorff topology) to some
subset $G$ of $U$ (cf Theorem 4.2 in [Mi]).
 By the main lemma, $G$ is isometric to $S^{d}(\frac{1}{2})$, and by
 Proposition \ref{properties of good points}, $G \subset G_U$.   In fact
 $G= G_U$ by the maximality of $d$.   So
 $L_{v}^{d}$ is closed.  To see that it is open, let $U$ in $L_{v}^d$ and
 let $W\in L_v$ be close
to $U$.  Consider the set $L(U, U')$ of lines between $U$ and $U'$.
Each $u\in G_U$ is on exactly one line $Z_u\in L(U,U')$, and the map
$G_U \longrightarrow G_W$ given by $u\longmapsto Z_u \cap W$ preserves
distances up to small additive error.  It follows
from the main lemma and the maximality of $d$
that if $W$ was originally chosen to be sufficiently close
to $U$, then $G_W$ is isometric to $S^{d}(\frac{1}{2})$.  $\square$

Consider the following
subset of $TM$.
\begin{eqnarray*}
TC|_V(\pi/2) =  \hspace{4in} \mbox{} \\
 \{ v\in TM|_{G_V} \ | \ \| v\| \leq \pi/2 \; \mbox{and $v$
		is tangent to a line in the configuration} \}.
\end{eqnarray*}
If the Radius Rigidity Theorem is false, then  $\exp :
TC|_V(\frac{\pi}{2})
\longrightarrow C$ is a surjective Lipshitz map, and the set of points in
M
whose inverse image is a singleton is an
open and  dense set.
 Indeed, $\exp$ is surjective
since the configuration has to cover $M$, and $\exp$ has
unique preimages on $M \backslash G_M$.
The set $M \backslash G_M$ is open and dense,
since $G_M$ consists of points of the form $U'$ where $U$ is a line
in the configuration, and the points of this form all lie in
proper subspheres of lines of type 2.

The fact that $\exp$ is surjective and Lipshitz yields a contradiction in
case $d \leq 3$
since it implies that
$dim_{Haus}\ M \leq dim_{Haus}\ TC|_V(\pi/2) \leq 3 + 3 + 8 = 14$.

The case $5\leq d \leq 7$ is also easy to eliminate since in this
 case $dim TC|_{V}(\frac{\pi}{2}) \geq 5 + 5 + 8 > 16 = dim\ M$.
So it is impossible for $\exp|_{TC|_{V}}$ to have unique preimages on an
open dense set.

The case $d=4$ is also not possible, but it is much harder to rule out.
We will get a contradiction in this case by showing that there is
 $(S^4 \times S^8)$-bundle $E$ over $S^4$ and a degree 1 map from $E$ to
 $M$.
To see that this is a contradiction, note that
a spectral sequence argument shows that if $E$ is any
$(S^{4} \times S^{8})$-bundle over $S^4$, then
\begin{eqnarray}\label{cohomology}
H^{*}(E) \cong H^{*}(S^4 \times S^{8} \times S^4).
\end{eqnarray}
  Since $H^*(M) \cong H^*(CaP^2)$, the existence of a degree $1$ map
$E \longrightarrow M$ implies that the fundamental cohomology
class in $E$ has a square root, and a simple algebraic argument
combined with (\ref{cohomology}) shows that it does not.

\begin{prop}\label{HP^2}
If $d=4$, then $G_M$ is a totally geodesic
submanifold of $M$ that
is isometric to $HP^2$ with its canonical metric with $1\leq sec\ HP^2
\leq 4$.
\end{prop}
{\em Proof:}
For any line $U$ in $C$ we can let $U$ play the role of $V$ and
define a configuration $C_U$ consisting of lines of type $1_U, 2_U$, and
$3_U$ in a way analogous to what we did on page \pageref{configuration}.
Of course assertion (\ref{C is everything}) is valid for each $C_U$, and
for each such configuration $C_U$,
$G_M =  \cup_{W \; \mbox{a line of type $2_U$ } } G_W$, since otherwise
there
would be a line not included in $C_U$.

Now let $u$ and $w$ be two points in $G_M$.  Since $w$ must lie on a line
of type $2_{ \{u\}'}$,
there is a Cayley line $Z$ containing $u$ and $w$.  Since $G_Z$ is
isometric to $S^4(\frac{1}{2})$, we can find a geodesic in $G_Z$ between
$u$ and $w$.  Using Lemma \ref{wiedersehen} and
the fact that $Z$ is totally $\pi$-convex
we see that if $dist(u,w) < \frac{\pi}{2}$ then the geodesic constructed
above
is the unique geodesic in $M$ between $u$ and $w$.  This shows that $G_M$
is
totally $\frac{\pi}{2}$-convex and hence, by $[CG1]$, a topological
manifold
with boundary and smooth, totally geodesic interior.  But the above
construction also indicates
that every geodesic in $G_M$ can be indefinitely prolonged (to a geodesic
in
$G_M$).  Therefore $\partial G_M = \emptyset$.  Thus $G_M$ with its
intrinsic
metric is a Riemannian manifold with sectional curvature $\geq 1$ and
diameter $= \frac{\pi}{2}$.  The proposition follows by analyzing the
structure
of the dual sets in $G_M$ and applying the classification theorem in
$[GG1]$.
$\square$

To construct $E$, first let $E' = \{ v \in TC|_{V}(\frac{\pi}{2}) \ | \
\| v \| = 1 \}$.  Let $p_{E'}: E' \longrightarrow G_V$ be the restriction
of the
projection map of $TC|_{V}(\frac{\pi}{2})$.  Given any $v \in G_V$,
$p_{E'}^{-1}(v) = \Pi_{v}^{-1}(G_{\{v\}'})$, where $\Pi_{v} : SM_{v}
\longrightarrow
\{v\}'$ is the Riemannian submersion discovered by Gromoll and Grove.
Since $G_{  \{ v\}' }$ is contractible in $\{ v\}'$,
$\Pi_{v}|_{ \Pi_{v}^{-1}(G_{\{v\}'}) }$ is trivial.  This shows that
that $E'$ is an $S^7 \times S^4$-bundle over $G_V$.  The desired bundle
$E$ will be obtained by suspending the ``$S^7$ parts'' of the fibers of
$E'$.  To help see that this can be done we prove
\begin{lem}\label{Q exisits}
There is a bundle $S^4 \hookrightarrow Q \stackrel{p_Q}{\longrightarrow}
G_V$ and a bundle $S^7  \hookrightarrow E'
\stackrel{p_{E',Q} }{\longrightarrow} Q$ so that
$p_{E'} = p_{Q} \circ p_{E',Q}$.
\end{lem}
{\em Proof:} Let $P = \{ v\in NG_V \subset TG_M \; | \; \| v\| \leq
\frac{\pi}{4} \} $, and let $Q$ be the double
of $P$ (cf [Mu]).
 For convenience we distinguish between the two copies
$P_1$ and $P_2$ of $P$ in $Q$ by setting
\begin{eqnarray*}
\| v\|_Q = \left\{
\begin{array}{ll} \hspace*{.14in}
\|v \|_{ P_1 } & \mbox{if $v \in P_1$} \\
-\| v\|_{ P_2 } +\frac{\pi}{2} &  \mbox{if $v\in P_2$}
\end{array}\right.
\end{eqnarray*}

For $i=1,2$, let  $p_{P_i} : P_i \longrightarrow G_V$ be the
the projection map of $P_i$.
By setting
\begin{eqnarray*}
p_Q (v) = \left\{ \begin{array}{ll}
   a \circ p_{P_1}(v) & \mbox{if $v\in P_1$} \\
a \circ p_{P_2}(v)    & \mbox{if $v\in P_2$}
\end{array} \right.
\end{eqnarray*}
 we see that $Q$ is an
 $S^4$-bundle over $G_V$.  (Here $a$ is the antipodal map of $G_V$.)

We can even define an exponential map $\exp_Q : Q \longrightarrow M$
by setting
\begin{eqnarray*}
\exp_Q v = \left\{
\begin{array}{ll}
\exp_{P}  \frac{ \| v\|_{Q} }{ \| v\|_{P} } v & \mbox{if
					   $\| v\|_Q \not=0,
					   \frac{\pi}{2}$ }\\
p_{P_1}(v) & \mbox{if $\| v\|_Q = 0$} \\
V'     & \mbox{if $\| v\|_Q = \frac{\pi}{2}$.}
\end{array} \right.
\end{eqnarray*}
Using the definition of $\| \cdot \|_Q$, Lemma \ref{wiedersehen},
and the definition of double ([Mu]), it is easy to check that $\exp_Q$ is
well defined even when $\| v\|_Q = \frac{\pi}{4}$ and
smooth even when $\| v\|_Q = 0, \frac{\pi}{4}, \frac{\pi}{2}$.

Let $p_{E', Q} : E' \longrightarrow Q$ be the map such that $p_{E',Q}(u) =
x$
if and only if $u\in T_vU$ for the Cayley line $U$ and $v\in G_V$, $x\in
p_{Q}^{-1}(v)$, and $\exp_Q x = U \cap \{ v\}'$.
$p_{E',Q}$ is smooth on $p_{E',Q}^{-1} ( \exp_{Q}^{-1}(G_M \backslash V')
)$
since on this set it is the composition of the smooth map
$\exp_{Q}^{-1} $ with the map $E' \longrightarrow M$ given by
$e \mapsto \exp \frac{\pi}{2} e$.  It is also clear that
\begin{description}
\item[(\romannumeral1)]
	$p_Q \circ p_{E',Q} = p_{E'}$, and
\item[(\romannumeral2)]
	the restriction of $p_{E',Q}$ to a fiber
  $p_{E'}^{-1}(v)$ of $p_{E'}^{-1}$ is $\exp_{Q}^{-1} \circ \Pi_{v} $.
\end{description}
Since $p_{E',Q}|_{p_{E',Q}^{-1} ( \exp_{Q}^{-1}(G_M \backslash V') ) }$
is smooth, it follows
from (\romannumeral1) and (\romannumeral2) that
it is a submersion.  Given any point
$x \in p_{Q}^{-1}(v) \subset Q$, $p_{E',Q}^{-1}(x) =
 E' \cap TU_v$ where $U$ is the Cayley line between
$\exp_Q x$ and $v$.  So $E'|_{p_{E',Q}^{-1} ( \exp_{Q}^{-1}(G_M \backslash
V') ) }$ is an $S^7$ bundle with projection map $p_{E',Q}$.

It remains to find
bundle charts for $p_{E',Q}$ about points in $\exp_{Q}^{-1}(V')$.
Let $\Phi$ be the geodesic flow of $M$.
Observe that
$$
p_{E',Q}^{-1} (\exp_{Q}^{-1}(V') ) = \{ u \in E' \ | \ u \; \mbox{is
tangent to
a line of type 2} \}
$$
 and that the map $p_{E',Q}^{-1} (\exp_{Q}^{-1}(V') ) \longrightarrow
 \Pi^{-1}_{V'}(G_V)$ given by
$u \mapsto \Phi^{\frac{\pi}{2}}(u)$ is a diffeomorphism which commutes
with the obvious projection
maps onto $G_V$.  This shows that $E'|_{p_{E',Q}^{-1} (\exp_{Q}^{-1}(V')
)}$ is a trivial $S^7$ bundle over $G_V$.  The $0$-section $s_0$ of $P_2$
provides a way of identifying
$\exp_{Q}^{-1}(V')$ with $G_V$, and $p_{E',Q}|_{p_{E',Q}^{-1}
(\exp_{Q}^{-1}(V') )} =
s_0 \circ a \circ p_{E'}|_{p_{E',Q}^{-1} (\exp_{Q}^{-1}(V') )}$; so there
is a diffeomorphism
$$
\psi: E'|_{p_{E',Q}^{-1} (\exp_{Q}^{-1}(V') )} \longrightarrow
\exp_{Q}^{-1}(V') \times S^7
$$
which commutes with the projections onto $\exp_{Q}^{-1}(V')$.

It follows from (\romannumeral2) that the restriction of $p_{E',Q}$ to any
fiber of $p_{E'}$ is a fiber bundle with fiber $S^7$.  In fact, given any
fixed $v\in G_V$, we can extend
$\psi|_{p_{E',Q}^{-1}(s_0 \circ a (v) ) }$ to a chart
\begin{eqnarray*}
p_{E',Q}^{-1} (P_2 \cap p_{Q}^{-1}(v)) \stackrel{\psi_v}{\longrightarrow}
(P_2 \cap p_{Q}^{-1}(v) ) \times S^7.
\end{eqnarray*}
by using the
holonomy dispacement maps of $\Pi_v$ given by the radial geodesics in
$\{ v\}'$ emanating from $V'$.  Clearly the $\psi_v$'s vary continuously
with $v$.  So given any open disk $U \subset G_V$ we get a bundle chart
\begin{eqnarray*}
p_{E',Q}^{-1} (P_2 \cap p_{Q}^{-1}(U) ) \stackrel{\psi_U}{\longrightarrow}
 P_2 \cap p_{Q}^{-1}(U) \times S^7
\end{eqnarray*}
by setting $\psi_U(x) = \psi_{p_{E'}(x)}(x)$.
 $\square$

Let $E$ be the $S^8$-bundle over $Q$ obtained by suspending
the fibers of $S^7 \hookrightarrow E' \stackrel{p_{E',Q}}{\longrightarrow}
Q$.
  We think of $E$ as the quotient space obtained
from $E' \times [0, \frac{\pi}{2}]$ by the equivalence relation
\begin{eqnarray*}
\begin{array}{ll}
(e_1, t) \sim (e_2, s) & \mbox{only if $t=s$.} \\
(e_1, 0) \sim (e_2, 0) & \mbox{if and only if $p_{E',Q}(e_1) =
p_{E',Q} (e_2)$} \\
(e_1, \frac{\pi}{2}) \sim (e_2, \frac{\pi}{2}) & \mbox{if and only if
$p_{E',Q}(e_1) = p_{E',Q} (e_2)$} \\
\mbox{for $0< t < \frac{\pi}{2}$},  (e_1, t) \sim (e_2,t) &
						\mbox{only if $e_1 = e_2$}
\end{array}
\end{eqnarray*}

Since $E$ is an $S^8$-bundle over $Q$ and $Q$ is an $S^4$-bundle over
$G_V$,
$E$ is an $S^4 \times S^8$-bundle over $G_V$.
We get the desired map $\psi :E \longrightarrow M$ by setting
$\psi[(e, t)] = \exp (t e)$.  (Here $[(e,t)]$ denotes the equivalence
class of $(e,t)$ in $E$.)  That $\psi$ is well defined
follows from Lemma \ref{wiedersehen}.  That $\psi$ is degree $1$ follows
from
the properties of $\exp: TC|_V(\frac{\pi}{2}) \longrightarrow M$.

\section{The Nonsimply Connected Case }
Let $M$ satisfy the hypotheses of Theorem \ref{pi_1}.
As we indicated in the introduction the classification theorem
in [GG1] applies to $M$.  In particular, we know that $M$
is either isometric to a space form or the quotient of $CP^{2k-1}$
in Theorem \ref{pi_1}(\romannumeral2).

 Suppose $M$ is a space form,
$O$ is an orbit of the action of $\Gamma$ on $S^{n}$, and
 $p: S^n \longrightarrow M$ is the universal covering map.  Then
since $rad\ M= \frac{\pi}{2}$, we can find a dual pair $(A, A')$ in
$M$ with $p(O) \in A$.  Since $A$ is totally $\pi$-convex so is
$\tilde{A} = p^{-1}(A)$.  Since $\partial A= \emptyset$ (2.5, 3.4,
and 3.5 in [GG1]), $\partial \tilde{A} = \emptyset$.  $\tilde{A}$
is therefore, a $\Gamma$-invariant, great subsphere of $S^n$ that
contains $O$.

On the other hand, if $S^n/ \Gamma$ is a space form and
an orbit $O$ of $\Gamma$ is contained
in a proper, invariant, totally geodesic subsphere, $S^k$, then
$A(S^k)= \{ x \in S^n \ | \ dist( x, S^k) = \frac{\pi}{2} \}$ is also
invariant.  Hence
$dist( p(O), p(A(S^k)) ) =\frac{\pi}{2}$, and $rad\  S^n / \Gamma =
\frac{\pi}{2}$, if every orbit of $\Gamma$ is contained in a proper,
invariant,
great subsphere.

To complete the proof of Theorem \ref{pi_1} it remains to show that
the space in (\ref{pi_1}, \romannumeral2) has radius $= \frac{\pi}{2}$.
The orbit of an arbitrary point for the
corresponding action on $S^{4k-1}$ is
$$
( z_1, z_2, \ldots , z_{2d} ) \mapsto ( \overline{z}_{d+1}, \ldots ,
\overline{z}_{2d}, - \overline{z}_{1}, \ldots , - \overline{z}_d ) \mapsto
(- z_1, \ldots , - z_d, -z_{d+1}, \ldots , - z_{2d} ) \mapsto
$$
$$
(- \overline{z}_{d+1}, \ldots , -\overline{z}_{2d},  \overline{z}_{1},
\ldots ,  \overline{z}_d )
\mapsto  (z_1, z_2, \ldots , z_{2d}).
$$
Thus each orbit (in $S^{4d -1}$) is contained in an invariant geodesic
that is
perpendicular to the fibers of the Hopf fibration $S^1 \hookrightarrow
S^{4d -1}
\longrightarrow  CP^{2d - 1}$.  It follows that each orbit in $CP^{2d -1}$
is contained
in an invariant geodesic.  If $\gamma : [0, \pi] \longrightarrow
CP^{2d-1}$ is an invariant
geodesic, then $(image(\gamma))'$ is also invariant (and $\not= \emptyset$
since $d\geq 2$).
So the radius of the quotient
is $\frac{\pi}{2}$.   $\square$

	Now we focus on the proof of (\ref{big space forms}).

\noindent
{\em Proof of (\romannumeral1):}
Let $\rho:\Gamma \longrightarrow O(n+1)$ be a fixed point free
representation
that respects an orthogonal splitting $V_1 \oplus V_2 \oplus \cdots \oplus
V_k$ of ${\Bbb R}^{n+1}$ so that $\rho|_{V_i}$ is irreducible for all i.
It follows from Theorem 7.2.18 in [W] that $dim\ V_i = dim\ V_j (\equiv
d)$
for all $i,j$.  Suppose $rad\ S^n/\rho(\Gamma) = \frac{\pi}{2}$.  Then
using (\ref{pi_1}, \romannumeral1) we can first
find a proper, invariant subspace that is not a direct sum of
$V_i$'s and then an irreducible invariant subspace $W$ that
is distinct from all of the $V_i$'s.  The orthogonal projections
$p_i :W \longrightarrow V_i$ are all $\rho$-equivariant, so
by Schur's Lemma,  they are either
zero maps or isomorphisms.  If they are all zero maps, then
we have $W \subset (V_1 \oplus V_2 \oplus \cdots \oplus V_k)^{\perp}$,
which is impossible.  So at least one of the projections (say $p_1$) is an
isomorphism.  If all of the other $p_i$'s are zero maps, then we have
$W \subset (V_2 \oplus  \cdots \oplus V_k)^{\perp}$, which would imply
that
$W= V_1$, also impossible.  So at least one other projection
(say $p_2$) is an isomorphism.
Thus $\rho|_W$ is linearly (and hence orthogonally
by Lemma 4.7.1 in [W]) equivalent to both $\rho|_{V_1}$ and
$\rho|_{V_2}$.

\noindent
{\em Proof of (\romannumeral2):}
By (\romannumeral4) it suffices to consider the case $k = 2$.
In this case the action of $\rho(\Gamma)$ on $S^3$ is orthogonally
equivalent to a subaction
$\bar{\rho}(\Gamma)$ of
the Hopf action, and the Hopf fibration $S^{1} \hookrightarrow S^3
\longrightarrow
S^2$ induces a Riemannian submersion $S^1 \hookrightarrow S^3 /
\bar{\rho}(\Gamma)
\longrightarrow S^2(\frac{1}{2})$.  Thus $rad\ S^3 / \rho(\Gamma) =
rad\ S^3 / \bar{\rho}(\Gamma) \geq  rad\ S^2(\frac{1}{2}) =
\frac{\pi}{2}$.

\noindent
{\em Proof of (\romannumeral3):}
Suppose there are points
$u, v \in S^{d-1}$ and members $g_1, g_2, \ldots,
g_{d}, g_{d + 1 }$ of $\Gamma$ so that
 $\{ g_1(u) , g_2(u), \ldots, g_{d}(u) \}$
and  $\{ g_1(v) , g_2(v), \ldots, g_{d}(v) \}$
are linearly independent and such that the
sets of coefficients
$\{ a_i\}$, $\{ b_i \}$ so that
\begin{eqnarray*}
a_1 g_1(u) + \cdots + a_{ d } g_{d}(u) =
g_{d +1 }(u)   \; \; \mbox{and}    \\
b_1 g_1(v) + \cdots + b_{d} g_{d}(v) =
g_{d + 1}(v)
\end{eqnarray*}
are distinct.
It follows that the vector $\frac{1}{\sqrt{2}} (g_{d+1}(u),
g_{d+1}(v))$ in ${\Bbb R}^{2d}$ is not in the span
of $\{ \frac{1}{\sqrt{2}}g_1( u, v), \frac{1}{\sqrt{2}}g_2 (u,v), \ldots
\frac{1}{\sqrt{2}}g_{d} (u,v) \}$. On the other hand,
 $\{ \frac{1}{\sqrt{2}}g_1( u, v), \frac{1}{\sqrt{2}}g_2 (u,v), \ldots
\frac{1}{\sqrt{2}}g_{d} (u,v) \}$ is linearly independent since its
projection onto ${\Bbb R}^{d} \times \{ 0\}$ is
 $\{ \frac{1}{\sqrt{2}}g_1 u, \frac{1}{\sqrt{2}} g_2 u, \ldots
\frac{1}{\sqrt{2}} g_{d} u \}$.  Therefore the set
 $$
\{ \frac{1}{\sqrt{2}}g_1( u, v), \frac{1}{\sqrt{2}}g_2 (u,v), \ldots
, \frac{1}{\sqrt{2}}g_{d} (u,v), \frac{1}{\sqrt{2}}g_{d+1} (u,v) \}
$$
is linearly independent, and we are be done by Theorem 4(\romannumeral1).
But according to [W], the image of every irreducible representation of a
fixed point free, nonabelian group $\Gamma$ contains matrices of the form
\begin{eqnarray}\label{type 1}
A =\left(
\begin{array}{cccccc}
e^{\frac{2\pi i k}{m} } & \mbox{}                       & \mbox{}  &
\mbox{}  & \mbox{}  & \mbox{}  \\
    \mbox{}             & e^{\frac{2\pi i kr}{m} }     & \mbox{}  &
    \mbox{}  & \mbox{}   & \mbox{}  \\
    \mbox{}             &               \mbox{}        & \cdot{}   &
    \mbox{}  & \mbox{}   & \mbox{}  \\
    \mbox{}             & \mbox{}                     & \mbox{}  &
    \cdot{}  & \mbox{}    & \mbox{}   \\
    \mbox{}             & \mbox{}                     & \mbox{}  &
    \mbox{}  & \cdot    & \mbox{}    \\
    \mbox{}             & \mbox{}                     & \mbox{}  &
    \mbox{}  & \mbox{}    &
e^{\frac{2\pi i kr^{d-1} }{m}  } \\
\end{array}\right)
\; \mbox{and}
\;
B =
\left(\begin{array}{cccccc}
0       &   1        &  \cdot    & \cdot     &  \cdot{} & 0     \\
\cdot   &   \mbox{}  &  \cdot    & \mbox{}   &  \mbox{}  & \mbox{}\\
\cdot   &   \mbox{}  &  \mbox{}  & \cdot     & \mbox{}   &\mbox{}\\
\cdot   &   \mbox{}  &  \mbox{}  & \mbox{}   & \cdot     & \mbox{} \\
   0    &   \mbox{}  &  \mbox{}  & \mbox{}   & \mbox{}    &   1     \\
e^{\frac{2\pi i l}{n'}}&0& \cdot  & \cdot     & \cdot      &  0  \\
\end{array} \right)
\end{eqnarray}
where $d, k, l, m, n'$ and $r$ are as in Theorem 5.5.6 in [W] (and we have
used
complex coordinates).  So it suffices
to set $u = (1, 0, \ldots , 0)$, $v = (0,1,0, \ldots , 0)$,
$g_1 = Id, g_2 =  A, g_3 = BA, g_4 = ABA,  \ldots , g_{d+1} =
(BA)^{\frac{d}{2} }$.
(Note that $d$ is even.)
 (To quickly see that there are matrices of the form (\ref{type 1})
in the image of every irreducible representation of a fixed point free,
nonabelian group,
note that such matrices are in the image of
every such representation of a so called ``group of type 1'' (Theorem
5.5.6 and
5.5.10 in [W]) and that
other nonabelian fixed point free groups contain groups of type 1
as subgroups ([W] pages 204-208).)     $\square$

\noindent
{\em Proof of (\romannumeral4):}
View $S^{n +d}(1)$ as the join $S^n(1) * S^{d-1}(1)$, and view
$\rho \oplus \sigma$ as the join of $\rho$ and $\sigma$.
Then every orbit of $\rho \oplus \sigma$ is contained in the join of an
orbit of $\rho$ with an orbit of $\sigma$.  Since the orbits of $\rho$ are
all
contained in proper invariant totally geodesic subspheres of $S^n$, the
orbits
of $\rho \oplus \sigma$ are contained in the joins of proper great
subspheres of $S^n$ with $S^{d-1}$.  $\square$

\noindent
{\em Proof of (\romannumeral5):}
For example, if $k$ is so large that the order of $\Gamma$ is less than
$n+1$, then every orbit is
automatically contained in an invariant subspace. $\square$ \vspace{.3in}


\noindent
{\Large {\bf Appendix} } \vspace{.2in}

Here we prove Lemma \ref{appendix lemma}.
Let $UVS^{15}$ and $UHS^{15}$ be the subbundles of $VS^{15}$ and
$HS^{15}$, respectively, that consist of unit vectors.
Set $X = \dot{ \tilde{\gamma} }$.
Let $\{ E_i\}_{i=1}^{7}$ and $\{ \tilde{E}_i\}_{i=1}^{7}$ be as in the
statement of Lemma
\ref{horizontal average lemma},
and define $A_{X}^{*}: VS^{15}_{ \tilde{\gamma}(t) } \longrightarrow
HS^{15}_{ \tilde{\gamma}(t) }$ by  $A_{X}^{*} v= A_X v$.
\newline
\begin{prop}\label{induction step}
Suppose $K_v \oplus L_v$ and $K_h \oplus L_h$ are orthogonal splittings of
$VS^{15}|_{\tilde{\gamma}(t) }$ and $HS^{15}|_{\tilde{\gamma}(t) }$
respectively so that $A_{X}^{*} K_v \subset K_h$, $A_{X}^*L_v \subset
L_h$,
$A_{X} K_h \subset K_v$, and $A_{X} L_h \subset L_v$.
Let $v\in L_v \cap UVS^{15}|_{ \tilde{ \gamma}(t) }$ be such that
$\| A_{X} v\| = \max_{ w\in L_v \cap UVS^{15}|_{\tilde{\gamma}(t)}}
\| A_X w\|$,
and suppose that $A_X v = \lambda y$, where $y$ is a
unit vector and $\lambda >0$.  Then
\begin{description}
\item(\romannumeral1)
	$A_X y = -\lambda v$,
\item(\romannumeral2)
	$\| A_X y\| =  \max_{z\in L_h\cap UHS^{15}|_{ \tilde{\gamma}(t) }
	}
\| A_X z\| $,
\item(\romannumeral3)
	$A_{X}^{*}(v^{\perp}) \subset y^{\perp}$, and
\item(\romannumeral4)
	$A_X (y^{\perp} \cap HS^{15}) \subset v^{\perp}$.
\end{description}
\end{prop}
{\em Proof:}
If (\romannumeral1) is false, then $A_X y = \lambda_w w$ for some
unit vector $w\in L_v$ that is different from $v$ and $ -v$.
Note that
$$
| \lambda_w | = | \lambda_w \langle w, w \rangle | >
|\lambda_w \langle w, v \rangle | =
| \langle A_X y, v \rangle |=
| \langle y, A_X v \rangle | = | \lambda |,
$$
and
$\langle A_X w, y \rangle =
- \langle w, A_X y \rangle = - \lambda_w$, which contradicts the
maximality
of $\| A_X v\|$.

If (\romannumeral2) is false, then there is a horizontal unit vector
$y_1 \in L_h$ different from $y$ and $-y$, a unit vector $w\in L_v$ and a
$\lambda_1 > \lambda$ so
that $A_X y_1 = \lambda_1 w$.  But then
$\langle A_X w, y_1 \rangle = - \langle w, A_X y_1\rangle = -\lambda_1$,
and this contradicts the maximality of $\|A_X v\|$.

If $w$ is a vertical vector in $v^{\perp}$, then
$\langle A_X w, y\rangle = - \langle w, A_X y \rangle = 0$, by part
(\romannumeral1).

Similarly, if $z$ is a horizontal vector in $y^{\perp}$, then
$\langle A_X z, v\rangle = - \langle z, A_X v\rangle = 0$.
\hspace*{5in} $\square$

Using Proposition \ref{induction step} we can inductively construct an
orthonormal basis for $VS^{15}_{\gamma(t)}$, $\{ u_1, \ldots, u_k,
u_{k+1},
\ldots, u_7\}$ with the following properties:
\begin{description}
\item[(\romannumeral1)]
$\{ u_1, \ldots, u_k\}$ is a basis for $ker\ A_{X}^{*}$,
\item[(\romannumeral2)]
$\{ u_{k+1}, \ldots, u_7\}$ satisfies $\langle A_X u_i, A_X u_j \rangle =
0$,
for $i\not= j$, and
\item[(\romannumeral3)]
for $i= k+1, \ldots 7 $, $A_X u_i = \lambda_i y_i$ for some unit vector
$y_i$ and some $\lambda_i > 0$, and $A_X y_i = - \lambda_i u_i$.
\end{description}

To construct $\{ u_{k+1}, \ldots , u_7 \}$ choose $u_{k+1}$ so that
$\| A_X u_{k+1} \| = \max_{w \in UVS^{15}_{\gamma(t)} } \| A_X w\| $
and in general choose $u_i$ so that $\| A_X u_i \| =
\max_{w\in (Span \{ u_1, \ldots , u_k, u_{k+1}, \ldots, u_{i-1}
\})^{\perp}
\cap UVS^{15}_{\gamma(t)} } \| A_X w\|$.

Let $\{ y_1, \ldots y_7\}$ be a completion of $\{ y_{k+1}, \ldots, y_7\}$
to an orthonormal basis for $HS^{15}_{\tilde{\gamma}(t)} \cap x^{\perp}$.
Note that
\addtocounter{dfn}{1}
\begin{eqnarray}\label{Ricci}
Ric (d\Pi (X), d\Pi (X) ) = \nonumber \\
 \Sigma_{i= 1}^{7} K(d \Pi (X), E_i ) =
\Sigma_{i=1}^{7} K(d\Pi(X), d\Pi(y_i))
\end{eqnarray}
Using the Horizontal Curvature Equation, (\ref{Ricci}) becomes
\addtocounter{dfn}{1}
\begin{eqnarray}\label{hor Ricci 1}
 \Sigma_{i= 1}^{7} K(X, \tilde{E}_i) + 3 \| A_X \tilde{E}_i \|^2 =
\Sigma_{i=1}^{7} K(X, y_i)  + 3 \| A_X y_i\| ^2.
\end{eqnarray}
And since the sectional curvature of $S^{15}$ is constant, (\ref{hor Ricci
1})
becomes
\addtocounter{dfn}{1}
\begin{eqnarray}\label{hor Ricci 2}
 \Sigma_{i= 1}^{7}  \| A_X \tilde{E}_i \|^2 =
\Sigma_{i=1}^{7} \| A_X y_i\| ^2.
\end{eqnarray}
Using (\ref{hor Ricci 2}) and the properties of the
$u_i$'s we see that
\addtocounter{dfn}{1}
\begin{eqnarray}\label{hor ricci}
\Sigma_{i =1}^{7} \| A_{X} \tilde{E}_i(t) \|^2 =
\Sigma_{i =1}^{7} \| A_{X} u_i \|^2.
\end{eqnarray}
For our given orthonormal basis $\{V_1(t), \ldots, V_7(t)\}$ for
$VS^{15}|_{\tilde{\gamma}(t)}$ we have
$A_X V_j(t) = \Sigma_i A_X \langle V_j(t), u_i \rangle u_i$.  Therefore
\begin{eqnarray*}
\Sigma_j \| A_{X} V_j(t) \|^2 =
\Sigma_j \Sigma_i \langle V_j(t), u_i \rangle^2 \|A_X  u_i\|^2= \\
\Sigma_i (\Sigma_j \langle V_j(t), u_i \rangle^2) \|A_X  u_i\|^2 =
\Sigma_{i =1}^{7} \| A_{X} u_i \|^2.
\end{eqnarray*}
Combining this with (\ref{average A-tensor bound}) and (\ref{hor ricci})
completes the proof of (\ref{appendix lemma}).

\vspace{.5in} \mbox{} \newline
\begin{center}
 Frederick Wilhelm \\
 Department of Mathematics    \\
  SUNY at Stony Brook        \\
 Stony Brook, NY 11794   \\
 email:  wilhelm@math.sunysb.edu
\end{center}


\begin{thebibliography}{Ber1}
\bibitem[Ber1]{Ber1} M. Berger, {\em Les vari\'{e}t\'{e}s riemanniennes
$\frac{1}{4}$-pinc\'{e}es},
Ann. Scuola Norm. Sup. Pisa {\bf 14}, (1960), 161-170.
\bibitem[Ber2]{Ber2} M. Berger, {\em Sur les vari\'{e}t\'{e}s \`{a}
courbure positive de
diam\`{e}tre minimum}, Comment. Math. Helv. {\bf 35}, (1961), 28-34.
\bibitem[Ber3]{Ber3} M. Berger, {\em Une borne inf\'{e}rieure pour le
volume d'une vari\'{e}t\'{e} riemannienne en fonction du rayon
d'injectivit\'{e}}, Ann. Inst. Fourier (Grenoble) {\bf 30}, (1980),
259-265.
\bibitem[Bw]{Bw} W. Browder, {\em Higher torsion in H-spaces}, Trans.
Amer. Math. Soc.  {\bf 108}, (1963), 353-375.
\bibitem[CG1]{CG1} J. Cheeger and D. Gromoll, {\em On the structure
of complete manifolds of nonnegative curvature}, Ann. of Math. {\bf 96},
(1972),
413-443.
\bibitem[CG2]{CG2} J. Cheeger and D. Gromoll, {\em On the lower bound for
the injectivity radius of $1/4$-pinched manifolds}, J. Differential
Geometry {\bf 15}, (1980), 437-442.
\bibitem[F]{F} T. Frankel, {\em Manifolds with positive curvature},
Pacific Journal of Mathematics {\bf 11}, (1961), 165-174.
\bibitem[GG1]{GG1} D. Gromoll and K. Grove, {\em A generalization of
Berger's rigidity
theorem for positively curved manifolds}, Ann. Sci. \'{E}cole Norm. Sup.,
4 s\'{e}rie, t. 20, 227-239.
\bibitem[GG2]{GG2} D. Gromoll and K. Grove, {\em The low-dimensional
metric foliations
of Euclidean spheres}, J. Differential Geometry {\bf 28}, (1988),
143-156.
\bibitem[GG3]{GG3} D. Gromoll and K. Grove, personal communication.
\bibitem[GP1]{GP1} K. Grove and P. Petersen, {\em On the
excess of metric spaces and manifolds}, preprint.
\bibitem[GP2]{GP2} K. Grove and P. Petersen, {\em Volume comparison
 \`{a} la Alexandrov}, Acta. Math., {\bf 169}, (1992), 131-151.
\bibitem[GS]{GS} K. Grove and K. Shiohama,
{\em A generalized sphere
theorem}, Ann. of Math. {\bf 106}, (1977), 201-211.
\bibitem[H]{H} R. Herman, {\em A sufficient condition that a mapping of
Riemannian manifolds be a submersion}, Proc. Amer. Math. Soc. {\bf 11},
(1960), 236-242.
\bibitem[K]{K} W. Klingenberg, {\em \"{U}ber Riemmannian
Mannigfeltigkeiten mit positiver Kr\"{u}mmung}, Comment. Math.
Helv. {\bf 35}, (1961), 47-54.
\bibitem[KS]{KS} W. Klingenberg and T. Sakai, {\em Injectivity radius for
$1/4$-pinched
manifolds}, Archiv der Math. {\bf 39}, (1980), 371-376.
\bibitem[Mi]{Mi} E. Michael, {\em Topologies on spaces of subsets}, Trans.
Amer.
Math. Soc. {\bf 71}, (1951), 152-182.
\bibitem[Mu]{Mu} J. Munkres, {\em Elementary Differential Topology},
Annals of Math. Studies no. 54, Princeton University Press, 1966.
\bibitem[O]{O} B. O'Neill, {\em The fundamental equations of a
submersion},
Michigan Math. J. {\bf 13}, (1966), 459-469.
\bibitem[Rch]{Rch} H. E. Rauch, {\em A contribution to Riemannian geometry
in the
large}, Annals of Math. {\bf 54}, (1951), 38-55.
\bibitem[Rj]{Rj} A. Ranjan, {\em Riemannian submersions of spheres with
totally
geodesic fibers}, Osaka J. of Math. {\bf 22}, (1985), 243-260.
\bibitem[Sa]{Sa} T. Sakai, {\em On the diameter of some Riemannian
manifolds},
Archiv der Math. {\bf 30}, (1978), 427-434.
\bibitem[Sh]{Sh} \mbox{K. Shiohama, {\em Topology of positively curved
manifolds with a certain }}
\newline
\mbox{ {\em diameter}, in \underline{Minimal Submanifolds and Geodesics:
Proceedings of the} }\newline
\underline{Japan-United States Seminar held in Tokoyo, September
19-22, 1977}, ed. M. Obata, North-Holland, 1977.
\bibitem[SS]{SS} T. Sakai and K. Shiohama, {\em On the structure of
positively
curved manifolds with certain diameter}, Math Z. {\bf 127}, (1972),
75-82.
\bibitem[SY]{SY} K. Shiohama and T. Yamaguchi, {\em Positively
curved manifolds with restricted diameters} in  \underline{Geometry of
Manifolds},
ed. K. Shiohama, Perspectives in Math. 8, Academic Press (1989), 345-350.
\bibitem[W]{W} J. Wolf, \underline{Spaces of Constant Curvature}, Publish
or
Perish, Wilmington, 1984.
\end{thebibliography}
\end{document}